\documentclass[11pt,reqno,tbtags]{amsart}
\usepackage[]{amsmath,amssymb,amsfonts,latexsym,amsthm,enumerate,color}
\usepackage{amssymb}
\usepackage[numeric,initials,nobysame]{amsrefs}

\renewcommand{\qedsymbol}{\ensuremath{\blacksquare}}
\numberwithin{equation}{section}

\newtheorem{theorem}{Theorem}[section]
\newtheorem*{theorem*}{Theorem}

\newtheorem{claim}[theorem]{Claim}

\newtheorem*{observation*}{Observation}
\newtheorem{corollary}[theorem]{Corollary}

\theoremstyle{definition}{

\newtheorem*{remark*}{Remark}
}

\newcommand{\R}{\mathbb R}

\newcommand{\N}{\mathbb N}

\newcommand{\one}{\boldsymbol{1}}
\newcommand{\eps}{\varepsilon}
\newcommand{\ra}{\rightarrow}
\newcommand{\ol}{\overline}

\newcommand{\E}{\mathbb{E}}
\renewcommand{\P}{\mathbb{P}}
\DeclareMathOperator{\Var}{Var}
\DeclareMathOperator{\Cov}{Cov}

\newcommand{\cF}{\mathcal{F}}

\newcommand{\dra}{\overset{d}{\ra}}

\renewcommand{\epsilon}{\varepsilon}

\newcommand{\Aut}{\operatorname{Aut}}

\date{}

\begin{document}
\title{On Fixation of Activated Random Walks}

\author{Gideon Amir }
\address{Gideon Amir\hfill\break Department of Mathematics, University of Toronto\\
Toronto ON, M5S 2E4, Canada}
\email{gidi.amir@gmail.com}
\urladdr{}

\author{Ori Gurel-Gurevich}
\address{Ori Gurel-Gurevich\hfill\break
Microsoft Research\\
One Microsoft Way\\
Redmond, WA 98052-6399, USA.}
\email{origurel@microsoft.com}
\urladdr{}

\begin{abstract}
We prove that for the Activated Random Walks model on transitive unimodular graphs, if there is fixation, then every particle eventually fixates, almost surely. We deduce that the critical density is at most 1.

Our methods apply for much more general processes on unimodular
graphs. Roughly put, our result apply whenever the path of each particle has an
automorphism invariant distribution and is independent of other particles' paths, and the
interaction between particles is automorphism invariant and local. This allows us to answer a
question of Rolla and Sidoravicius \cite{Rolla, RS}, in a more
general setting then had been previously known (by Shellef \cite{Shellef}).
\end{abstract}

\maketitle

\section{Introduction}

Let $G=(V,E)$ be an infinite, locally finite transitive graph. Let $\mu$ be a distribution on $\N$. The \emph{activated random walk} (ARW) model is an interacting particle system on $G$, in which the number of particles initially on each vertex is independently distributed $\mu$. At the onset, all the particles are \emph{active}. Each active particle independently performs a continuous time rate 1 simple random walk on $G$ and also becomes \emph{inactive} independently at rate $\lambda$. An inactive particle stays put as long as there are no other particles at the same vertex. Once there is another particle (necessarily active) at the vertex of an inactive particle, the inactive particle immediately becomes active. In particular, as long as there are two (or more) particles at a vertex, none of them become inactive.

We say that a vertex $v$ \emph{fixates} at time $t$, if $t$ is the minimal time such that there are no active particles at $v$ after $t$. If $v$ fixates at some time $t$, we say it is \emph{fixating}, otherwise it is \emph{nonfixating}. Similarly, define when a particle fixates at time $t$, in which case it is fixating and otherwise it is nonfixating. We say that the model is \emph{vertex fixating} if every vertex fixates a.s.\  and \emph{particle fixating} if every particle fixates almost surely.

%If some vertex is nonfixating then the adjacent vertices are also
%nonfixating (with probability $1$), so \emph{vertex-fixation}, that
%is, the property that every vertex is fixating, is a global
%property. \ori{true in what generality? ignore this?}\gidi{Well its
%not alwys true since when the walks are not aperiodic things might
%go wrong, but I think we can ignore this} By ergodicity of the
%process, the probability of vertex-fixation is either 0 or 1.

%\ori{is particle fixation global? yes, in
%our case, because non-particle-fixating implies non-vertex-fixating
%implies no particle fixates. But not for more general interactions,
%e.g. if a particle falls asleep when alone and never wakes
%up}\gidi{ignore this}

Recall that a transitive graph is \emph{unimodular}, if for any two
vertices $x$,$y$ we have $|A(x,y)|=|A(y,x)|$, where $A(x,y)=\{f(y)\ |\
f\in \Aut(G) , f(x)=x\}$. These include all amenable transitive
graphs, all Cayley graphs as well as all regular trees. For
unimodular graphs the following symmetry principle holds:

%\ori{elaborate}\gidi{I am not sure we need to elaborate more}.

\textbf{Mass Transport Principle} :\  Let $G = (V, E)$ be a
unimodular graph, and assume $F : V \times V \rightarrow [0,1]$ is
automorphism invariant (i.e. $F(\gamma v, \gamma y)=f(x,y)$ for any
$\gamma\in \Aut(G)$), then for any $v\in V$ we have
$$\sum_{w\in V} F(v,w) = \sum_{w\in V} F(w,v)$$
\\

Loosely speaking, the mass transport principle says that under
enough symmetry, the amount of mass transmitted from a vertex $v$ is
equal to the amount of mass entering $v$. For proofs, applications
and generalizations to other graphs, the reader is referred
to chapter $8$ of \cite{LP}.

Before stating our results, we wish to generalize the setting to include other behaviors of the particles.
A \emph{path} is a function $y:\R^+ \ra V$ such that there is an infinite sequence
$0=t_0<t_1<\ldots \ra \infty$ and $v_i \in V$ satisfying $y|_{[t_i,t_{i+1})} \equiv v_i$ and $(v_i,v_{i+1}) \in E$.
This is the path of the particle parameterized by its \emph{inner time}.
A distribution $Y$ on paths beginning at some fixed vertex $y(0)$, is \emph{invariant} if for every
$\gamma \in \Aut(G,y(0))$ we have $\gamma \circ Y \simeq Y$.
A distribution $Y$ on paths has \emph{infinite range} if the range of a path sampled from this distribution is a.s. infinite.
Obviously, the distribution of a continuous time, simple random walk has infinite range.

%
%
%A distribution $Y$ on paths is \emph{unpredictable} if for every $T\in \R^+$ there is a
%function $p^*_T:\R^+ \ra [0,1]$, %\gidi{I think p is not a good name, at least add a star}
%such that $p^*_T(t) \underset{t\ra \infty}{\longrightarrow} 0$ and for every $y:[0,T] \ra V$ we have $\max_{v\in V} \P(Y(t)=v \mid Y|_{[0,T]}=y)
%\le p^*_T(t)$. That is, knowing how
%the path started, one cannot give a good prediction for where it'll
%be in the future. Obviously, the path of a continuous time, simple
%random walk is invariant and unpredictable.

Given the path distribution, one can define the corresponding
\emph{generalized activated random walk} model by randomizing a
Poisson process indicating at which times the particle becomes
inactive if it is alone on a vertex. Similarly, one may define other
models by using different interaction rules. We may consider any
interaction which affects only the rate in which particles move
along their paths (including rendering particles inactive). Thus, an
\emph{interaction rule} for a particle is a function from all
possible states of the system at time $t$ (that is, the actual paths
of all particles until time $t$) into $\R^+$. The output of this
function determines the rate in which the actual path of the
particle progresses along its putative path. In order to account for
random interactions, we equip each particle with an independent
$U[0,1]$ random variable, which the interaction rule is also given
as input.

In this paper we will only be interested in automorphism invariant interaction rules, which are also \emph{local}, that is, the function actually depends only on those particles at some finite radius around our particle.

Of course, for some path distributions and interaction rules the
resulting model might not be well defined or unique or local (i.e.
that every vertex is visited by only finitely many particles in any
finite time interval). We do not concern ourselves with these
questions in this paper. Rather, we assume that the model is well
defined, unique and local.

\begin{theorem} \label{fixation}
If an invariant, infinite-range, local model on a unimodular graph is vertex-fixating, then it is particle-fixating.
\end{theorem}

Using Theorem~\ref{fixation} we immediately get the the following corollary,
which answers a question of Rolla and Sidoravicius \cite{Rolla, RS}.

\begin{corollary} \label{density}
If an invariant, infinite-range generalized ARW on a unimodular graph is
vertex-fixating, then $\E(\mu)\le 1$.
\end{corollary}
\begin{proof}
Let $F(u,v)$ be the probability that some particle, starting at $u$ fixates at $v$.
Obviously, $F$ is automorphism invariant.
Since only 1 particle may fixate at a given vertex we have $\sum_{u\in V} F(u,v) \le 1$.
By the mass transport principle this is equal to $\sum_{v\in V} F(u,v)$ which is therefore also at most 1.
By Theorem~\ref{fixation}, all particles fixate, so this sum is equal to $\E(\mu)$ which is therefore at most 1.
\end{proof}

Of course, the corollary applies for any model in which only 1 particle may fixate at any given vertex.

For the original ARW model, corollary \ref{density} was proved independently by Shellef \cite{Shellef}, using completely different methods, on any bounded degree graph. The main merit of our proof is that although the graph has be unimodular, it works for a broad class of random walks and interaction rules.

%This corollary was already proved, using completely different
%methods, by Shellef \cite{Shellef}\gidi{for what graphs?}. His proof
%uses some properties that are rather specific to ARW, namely,
%\emph{monotonicity} and \emph{commutativity}. \ori{elaborate}
%\gidi{Actually, the new proof doesn't really uses only the comparison
%with MSIA, and uses these properties in MSIA}
%
%The main merit of our proofs is that they work for much more general processes then the activated random walk model.
%Our only requirement is that each particle's movement is an independent simple random walk,
%so that the interaction between particles determines when they get activated and deactivated.
%\ori{actually, it's more general time change, and the SRW can also be relaxed}.
%In particular, we do not require neither monotonicity nor commutativity.
%
%\ori{open questions?}\gidi{saying something on density less than or
%equal to one?, i think we can skip}

\section{Proof of Theorem~\ref{fixation} }

We will assume the system is not particle-fixating, and show that the system cannot be vertex-fixating.

Let $N_v$ be the number of particles starting at $v$, which is
distributed $\mu$ independently for each vertex. We number the
particles and name the $i$-th particle starting at vertex $v$ by
$x^{v,i}$. For each of these particles a \emph{putative path}
$y^{v,i}_t$ is independently randomized using $\nu$, the invariant,
 infinite-range distribution on paths. The
real path of the particle $x^{v,i}_t$ follows the putative path with
a time change determined by the interaction with other particles.

Let $A^{v,i}$ be the event ``$N_v>i$ and $x^{v,i}$ is nonfixating''. Since the system is not particle-fixating, for some
$i$ we have $\P(A^{v,i})>0$. In fact, we have $\P(A^{(v,0)})>0$, since conditioned on $N_v>i$,
particles $x^{v,i}$ and $x^{v,0}$ are exchangeable, that is, they have the same conditional distribution,
and hence the same conditional probability of not fixating. Let $a=\P(A^{v,0})$.

Now, since our system is local, $A^{v,0}$ can be $\eps$-approximated by another event $B^{v,0}$ which depends only on the behavior of those particles starting in some finite radius $R$ around $v$ and up to some finite time $T$. If $B^{v,0}$ occurs we call $x^{v,0}$ a \emph{candidate}. A candidate which is actually non-fixating is \emph{good} and the rest of the candidates are \emph{bad}. By our approximations, the probability that $x^{v,0}$ is a bad candidate is bounded by $\eps$. Let $b=\P(B^{v,0}) \ge a-\eps$ which is positive for $\eps$ small enough.

Let $\cF_T$ be the
$\sigma$-algebra consisting of $\{N_v\}_{v\in V}$ and
$\{y^{v,i}|_{[0,T]}\}_{v\in V, i\in \N}$. This $\sigma$-algebra
determines the process of candidates, which is a $2R$-dependent
process since $B^{v,0}$ depends only on $y^{u,i}|_{[0,T]}$ for $u$
of distance at most $R$ from $v$. Fix some integer $n$ and for each
particle let $z^{v,0}$ be a vertex chosen uniformly and
independently from the first $n$ distinct vertices in the path $y^{v,0}$ (there are $n$ distinct vertices in the path a.s. since the path distribution is infinite-range).

Let $C(v,u)$ be the event ``$B^{v,0}$ and $z^{v,0}=u$''. Obviously, for any $v$ and $u$, the probability of $C(v,u)$ is at most $1/n$, even conditioned on $\cF_T$. Let
Let $q(v,u)=\P(C(v,u)|\cF_T)$ and let $Q(u)=\sum_{v \in V} q(v,u)$.

\begin{claim} \label{expectation}
$Q(u) \dra b$ as $n \ra \infty$.
\end{claim}
\begin{proof}
Consider the mass transport $F:V \times V \ra [0,1]$ defined by $F(v,u)=\P(C(v,u))$. It is automorphism invariant, so by the mass transport principle,
$$\E(Q(u))=\sum_{v \in V} F(v,u) = \sum_{u \in V} F(v,u) = \P(B^{v,0}) = b \ .$$

Since
$$\Var(q(v,u)) \le \E((q(v,u))^2) \le \E(q(v,u))/n$$
we have
$$\sum_{v \in V} \Var(q(v,u)) \le b/n \ra 0 \ .$$

Similarly, $\Cov(q(v,u), q(w,u)) \le \E(q(v,u))/n$ if $v$ and $w$ are at most $2R$ apart, and $\Cov(q(v,u), q(w,u))=0$ otherwise. Hence,
$$\sum_{v \in V} \sum_{w \in V} \Cov(q(v,u), q(w,u)) \le \sum_{v \in V} d^{2R} \frac{E(q(v,u))}{n} \le \frac{d^{2R}b}{n} \ra 0$$
where $d$ is the degree of a vertex in $G$. Put together, we get $\Var(Q(u))\ra 0$, completing the proof.
\end{proof}

Let $Z(u)=\sum_{v \in V} \one_{C(v,u)}$, that is, the number of
candidates for which $z^{v,0}=u$.

\begin{claim} \label{poisson}
$\lim_{n\ra\infty} \P(Z(u)>0 |\cF_T) \ge 1-e^{-b} \ .$
\end{claim}
\begin{proof}
$\P(Z(u)=0 |\cF_T)=\prod_{v \in V} (1-q(v,u)) \le e^{-Q(u)} \ra e^{-b}$
\end{proof}

Let $D(v,u)$ be the event ``$B^{v,0}$ and not $A^{v,0}$ and $z^{v,0}=u$'', and let
$\ol{Z}(u)=\sum_{v \in V} \one_{D(v,u)}$, that is, we count the number of bad candidates for which $z^{v,0}=u$.

\begin{claim} \label{bad}
$\E(\ol{Z}_t(u)) \le \eps$
\end{claim}
\begin{proof}
Consider the mass transport $F:V \times V \ra [0,1]$ defined by $F(v,u)=\P(D(v,u))$. It is automorphism invariant, so by the mass transport principle,
$$\E(\ol{Z}(u))=\sum_{v \in V} F(v,u) = \sum_{u \in V} F(v,u) = \P(B^{v,0} \cap \ol{A^{v,0}}) \le \eps$$
where the last inequality is due to the fact that bad candidates are rare.
\end{proof}

It follows from lemma \ref{poisson} and \ref{bad} that the
probability of the event ``There exists a good candidate $x^{v,0}$
such that $z^{v,0}=u$'' is at least $\P(Z(u)>0)-\P(\ol{Z}(u)>0)\ge
1-e^{-a+\eps}-2\eps$ when $n$ is large enough, which is positive
when $\eps$ is small enough. But if $z^{v,0}=u$ and $x^{v,0}$ is
nonfixating, then there is some $t$ such that $x^{v,0}_t=u$. In
particular, if the above event happens then the vertex $u$ has not
fixated until time $t$.
For any $T'$, the probability that $t<T'$ tends to $0$ as
$n\ra\infty$. Hence, the probability of any vertex to fixate before
time $T'$ does not tend to 1. Since the above analysis holds
conditioned on any prefix of the system, the probability of a vertex
to fixate is 0. \qedsymbol

\end{document}